\documentclass[a4paper]{article}

\usepackage{amsmath}
\usepackage{amssymb}
\usepackage{latexsym}\usepackage[all]{xy}
\usepackage{makeidx}
\makeindex
\newtheorem{Theorem}{Theorem}%[subsection]
\newtheorem{Lemma&Def}{Lemma and Definition}[Theorem]

\newtheorem{Lemma}[Theorem]{Lemma}

\begin{document}

\title{On accessibility of finitely generated groups}

\author{Richard Weidmann}

\date{}

\maketitle

\maketitle
\begin{abstract} We prove an accessibility result for finitely generated groups that combines Sela's acylindrical accessibility with Linell accessibility.
\end{abstract}
\section*{Introduction}

Grushko's theorem \cite{G} says that for any free product $G=G_1*\ldots *G_k$ we have the equality
 $$\hbox{rank }G=\sum\limits_{i=1}^k\hbox{rank }G_i$$  where the rank is the minimal number of elements needed to generate $G$. This implies in particular that any free product decomposition of a $k$-generated group $G$ has at most $k$ factors. Differently said this means that for any graph of group decomposition $\mathbb A$  of $G$ with non-trivial vertex groups and trivial edge groups we have $$b(A)+\#VA\le \hbox{rank }G$$ where $b(A)$ is the Betti number of the underlying graph $A$ and $\#VA$ is the cardinality of the vertex set $VA$.
 
\medskip Thus the rank of a group $G$ bounds the complexity of splittings of $G$ with trivial edge groups. Whether such bounds on the complexity exist under less restrictive assumptions has been a matter of much attention, these problems are generally called accessibility questions. It is clear that without making strong assumptions on the splittings no positive result can be true. 

\medskip Recalling that Grushko's theorem deals with actions on trees that have trivial edge stabilizers suggests two ways to relax these conditions:

\begin{enumerate}
\item Allow non-trivial edge stabilizers that still satisfy some assumptions.
\item Allow arbitrary edge stabilizers but demand stabilizers of long segments to be trivial.
\end{enumerate}

Both of these situations have been dealt with successfully:

\medskip
\noindent 1) The most basic result is due to P.~Linell \cite{L} who shows that the complexity of a splitting of a finitely generated group is bounded provided that all edge groups are finite with a uniform bound on the order.
In the case of finitely presented groups it has been shown by M.~Dunwoody \cite{Du1} that such a bound exists provided the splitting is reduced and all edge groups are finite. This has then been generalized to small splittings by M.~Bestvina and M.~Feighn \cite{BF1}. The assumption that the groups are finitely presented cannot be dropped, there are finitely generated groups without this property, see \cite{BF2} and \cite{Du2}.

\medskip\noindent 2) The second situation has been dealt with by Z.~Sela who proved the existence of a bound on the complexity of  a splitting of a finitely generated group provided the stabilizer of any segment of length $k$ is trivial for some $k$. This result is for example relevant in the construction of the JSJ-decomposition of limit groups where the group is not known to be finitely presented and the edge groups are infinite. An explicit bound of the complexity in terms of $k$ and the rank of $G$ was given in \cite{W3}.%The above mentioned accessibility results  don't help in this situation. %However Z.~Sela's accessibility result which bounds the complexity of acylindrical splittings does apply in this situation, here a splitting is called $k$-acylindrical if the stabilizer of any segment of length $k$ in the Bass-Serre tree is trivial. 

%\medskip Another accessibility result for finitely generated groups is that of P.~Linell that provides a bound on the complexity of  splittings eith the property that the edge groups are finite and that there exists a uniform bound on the order. In this note we prove a result that combines Linell accessibility and Sela's acylindrical accessibility. We further see that the acylindricity assumption can be projectivized.

\medskip In this note we prove an accessibility result  that is to Linell's theorem what Sela's result is to Grushko's theorem. 

\medskip Let $G$ be a group acting on a simplicial tree $T$ and $$\pi:T\to G\backslash T$$ be the canonical projection. Let $[v,w]$ be a simplicial segment of $T$. We say that $[v,w]\subset T$ has projective length $k$ and write $pl([v,w])=k$ if $\pi([v,w])$ meets $k$ edges (edge pairs) of $G\backslash T$.

\medskip Let $\mathbb A$ be a graph of groups, $G=\pi_1(\mathbb A)$ and $T$ be the associated Bass-Serre tree. We say that $T$ or $\mathbb A$ is $(k,C)$-acylindrical if the stabilizer of any segment $[v,w]\subset T$ with $pl([v,w])> k$ is of order at most $C$. We further call a graph of groups weakly reduced iff there exists no valence 2 vertex $v$ such that both boundary monomorphisms into the vertex group $A_v$ are isomorphisms.

\medskip We prove the following theorem, here $\#EA$ is the number of edges of the graph $A$ underlying the graph of groups $\mathbb A$.

\begin{Theorem}\label{main} Let $\mathbb A$ be a weakly reduced minimal  $(k,C)$-acylindrical graph of groups with $k\ge 1$. Then $$\#EA\le (2k+1)\cdot C\cdot (\hbox{rank }\pi_1(\mathbb A)-1).$$
\end{Theorem}

%Note that apart from being a tool to generalize various arguments from the torsion-free case to the case of groups with bounded torsion the projective component of this Theorem has the added benefit that

Note that the theorem does not hold for $k=0$. To see this note that the free group of rank $n$ is the fundamental group of a weakly reduced minimal graph of trivial groups with $3n-3$ edges, this splitting is clearly $(0,1)$-acylindrical. However, as any $(k,C)$-acylindrical splitting is also a $(k+1,C)$-splitting our theorem does still give a bound for the case $k=0$.

\medskip The proof of the theorem combines ideas of the author's proof of acylindrical accessibility \cite{W3} and M.~Dunwoody's proof of Linell accessibility \cite{Du3} in the language of \cite{W2}. It applies Stallings folding arguments and is not very hard although the complexity used in the argument is somewhat involved.

\medskip There are two other proofs that one could attempt to generalize to account for the present theorem:

\medskip The first one is Z.~Sela's original proof which uses the Rips machine \cite{S} and is therefore much less elementary. Although this might be possible it is not clear how this could account for the projective component of our statement and it would also not yields an explicit bound on the complexity but only the existence of a bound depending on $G$.

\medskip The second proof is that of T.~Delzant who first proves acylindrical accessiblity using his version of the Dunwoody resolution Lemma \cite{De} and then using the fact that any finitely generated freely indecomposable group $G$ has a finitely presented cover $H$ and an epimorphism $\phi:H\to G$ that does not factor through a free product \cite{Sw}. It is not clear to the author how to adapt this argument to the present case. Even if possible it would probably not yield an explicit bound on the complexity.

\medskip It seems that the ideas presented in this note together with the proof of BF-accessibility should also yield an acylindrical version of accessibility for finitely presented groups where the condition is that stabilizers of long segments are small. Indeed Dunwoody's proof of Linell accessibility and the proof of BF-accessibility are very similar in that they start a folding sequence with a well-understood splitting and gradually make it more complicated, on the way the complexity of the edge groups is being measured. In the case of Linell accessibility the original splitting is a wedge of circles and the measure of complexity for the edge groups is their order while in the proof of BF-accessibility the original splitting is a Dunwoody resolution and the measure of complexity is the complexity of the action on the Bass-Serre tree of the Dunwoody resolution. 

\medskip In Section~\ref{foldings} we briefly recall the folding machinery before we give the proof of Theorem~\ref{main} in Section~\ref{proof}.

\section{Foldings}\label{foldings}

In this section we briefly fix notations for graphs of groups and describe how a graph of groups can be approximated by a sequence of graphs of groups that are related by folds. We use the language of \cite{W4}, \cite{KMW} which essentially is an alternative formulation of the combination of foldings and vertex morphisms discussed in \cite{Du3} where M.~Dunwoody builds on earlier work of J.~Stallings \cite{St}, \cite{St2} and M.~Bestvina and M.~Feighn \cite{BF1}. We only recall some aspects of the theory, for more details we refer the reader to \cite{W4} and \cite{KMW}.

\medskip A {\em graph} $A$ consists of a set of vertices $VA$, a set of edges $EA$, an inversion ${^{-1}}:EA\to EA$ and maps $\alpha:EA\to VA$ and $\omega:EA\to VA$ such that $\alpha(e)=\omega({e^{-1}})$ for all $e\in EA$. We denote the Betti number of $A$ by $b(A)$. Recall that the Betti number of $A$ is the number of edge pairs outside a maximal subtree of $A$.

\medskip A {\em graph of groups} $\mathbb A$ consists of a graph $A$, vertex groups $A_v$ for every $v\in VA$, edge groups $A_e$ for every edge $e\in EA$ such that $A_e=A_{e^{-1}}$ and boundary monomorphisms $\alpha_e:A_e\to A_{\alpha(e)}$ and $\omega_e:A_e\to A_{\omega(e)}$ satisfying $\alpha_e=\omega_{e^{-1}}$; thus only the maps $\alpha_e$ need to be specified.

\smallskip Let $\mathbb A$ be a graph of groups. An $\mathbb A${\em -graph}
  $\mathcal B$ consists of an underlying graph $B$ with the following
  additional data:

\begin{enumerate}

\item A graph-morphism $[\,.\, ]: B\to A$.

\item For each $u\in VB$ there is a group $B_u$ with
  $B_u\le A_{[u]}$.
  
\item For each $f\in EB$ there is a group $B_f$ with
  $B_f=B_{f^{-1}}\le A_{[f]}$.
  
\item To each edge $f\in EB$ there are two associated group elements
  $f_{\alpha}\in A_{[\alpha (f)]}$ and $f_{\omega}\in A_{[\omega(f)]}$ such that
  $( f^{-1})_{\alpha}=(f_{\omega})^{-1}$ for all $f\in EB$.
  
\item For each $f\in EB$ we have $f_{\alpha}\cdot \alpha_{[f]}(B_f)\cdot f_{\alpha}^{-1}\le B_{\alpha(f)}$.
\end{enumerate}

When representing $\mathbb A$-graphs by figures we will label the vertices and edges as follows: If $f\in EB$ and $u\in VB$, we shall refer to  $(B_u,[u])$ as the label of $u$ and to $(f_{\alpha},[f], f_{\omega})$ as the label of $f$. Note that the graph with its  labels completely determine the $\mathbb A$-graph except the edge groups.%Note that to give the complete information we would need to speak of edge labels of type $(f_{\alpha},[f],B_f, f_{\omega})$ but we will usually not do this.

To any $\mathbb A$-graph we can then associate in a natural way a
graph of groups.
  Let $\mathcal B$ be an $\mathbb A$-graph. The {\em associated graph of
  groups} ${\mathbb B}$ is defined as follows:

\begin{enumerate}
\item The underlying graph of $\mathbb B$ is the graph $B$.

\item The vertex and edge groups are the groups $B_u$ for $u\in VB$ and $B_f$ for $f\in EB$.

\item For each $f\in EB$ we define the boundary monomorphism $\alpha_f:B_f\to
  B_{\alpha(f)}$ as
  $\alpha_f(g)=f_{\alpha}\big(\alpha_{[f]}(g)\big)f_{\alpha}^{-1}$ and $\omega_f=\alpha_{ f^{-1}}$.
\end{enumerate}

\smallskip 
For any $u_0\in VB$ and $v_0=[u_0]$ the $\mathbb A$-graph determines a group homomorphism $\nu:\pi_1(\mathbb B,u_0)\to \pi_1(\mathbb A,v_0)$ given by
%Let $\mathcal B$ be an $\mathbb A$-graph defining a graph of groups
%$\mathbb B$.
%Suppose $u,u'\in VB$ and $p$ is a $\mathbb
%  B$-path from $u$ to $u'$. Thus $p$ has the form:
\[
[b_0,f_1,b_1,\dots, f_s,b_s]\mapsto [(b_0g_1), e_1, (k_1b_1g_2), e_2,\dots, (k_{s-1}b_{s-1}g_s),e_s,
(k_sb_s)].
\]

We call $\nu(\pi_1(\mathbb B))$ the subgroup of $\pi_1(\mathbb A)$ the subgroup represented by $\mathcal B$ and say that $\mathbb A$-graph is surjective if $\nu$ is surjective. To any generating tuple $S$ of $\pi_1(\mathbb A)$ we can associate a surjective $\mathbb A$-graph called the $S$-wedge. The underlying graph of an $S$-wedge is a wedge of circles, at most one vertex group is non-trivial and all edge groups are trivial, in particular every hyperbolic element of $S$ corresponds to one of the circles.

\medskip There is further the notion of a folded $\mathbb A$-graph. For a folded $\mathbb A$-graph $\mathcal B$ the map $$b_0,f_1,b_1,\dots, f_s,b_s\mapsto (b_0g_1), e_1, (k_1b_1g_2), e_2,\dots, (k_{s-1}b_{s-1}g_s),e_s,
(k_sb_s)$$ maps reduced $\mathbb B$-paths to reduced $\mathbb A$-paths, in particular $\nu$ is injective. If $\mathcal B$ is both folded and surjective then it follows that $\mathbb B$ is isomorphic to $\mathbb A$. If an $\mathbb A$-graph is not folded then a fold can be applied to yield a new $\mathbb A$-graph. Often one can guarantee that finitely many folds transform a given $\mathbb A$-graph into a folded one. However this is not true in the context of this paper.

\medskip Following \cite{BF1} there are six types of folds, three of A-type and three of B-type. We will only discuss the A-type folds as the arguments for A-type folds and B-type folds are very similar and usually one can even restrict attention to A-type folds. After possibly applying auxiliary moves first (see \cite{KMW}) the three types of folds have the following effects on a $\mathbb A$-graph. For us important is that auxiliary moves preserve the isomorphism class of the underlying graph of groups. Note that any two $\mathbb A$-graphs that are related by a fold represent the same subgroup of $\pi_1(\mathbb A)$.

\medskip\noindent {\bf Fold of type IA: }
  Suppose $f_1$ and $f_2$ are two distinct non-loop edges and that
  $y=\omega(f_1)\ne \omega(f_2)=z$. Furthermore $a=(f_1)_{\alpha}=(f_2)_{\alpha}$ and $b=(f_1)_{\omega}=(f_2)_{\omega}$.

  We identify the edges edges $f_1$ and $f_2$ of $\mathcal B'$ into a
  single edge $f$ with edge group $\langle B_{f_1},B_{f_2}\rangle$ and label $(a,e,b)$. We set the label of $\omega(f)$
  to be
\[
(\langle B_y, B_z\rangle, v).\]

We call this operation a {\em fold of type}
IA.

\begin{figure}[here]
\centerline{ \footnotesize \setlength{\unitlength}{.9cm}
\begin{picture}(10,3)
\put(0,1.5){\line(3,1){3}}
\put(1.5,2){\vector(3,1){.01}}
\put(1.5,1){\vector(3,-1){.01}}
\put(0,1.5){\line(3,-1){3}}
\put(4.5,1.5){\vector(1,0){1.5}}
\put(5,1.65){IA}
\put(7.5,1.5){\line(3,0){3}}
\put(9,1.5){\vector(1,0){.01}}
\put(0,1.5){\circle*{.07}}
\put(3,2.5){\circle*{.07}}
\put(3,0.5){\circle*{.07}}
\put(7.5,1.5){\circle*{.07}}
\put(10.5,1.5){\circle*{.07}}
\put(-1.3,1.4){$(B_x,w)$}
\put(3.15,2.4){$(B_z,v)$}
\put(3.15,0.4){$(B_y,v)$}
\put(1,2.3){$(a,e,b)$}
\put(1,0.6){$(a,e,b)$}
\put(6.5,1.65){$(B_x,w)$}
\put(10,1.65){$(\langle B_y,B_z\rangle,v)$}
\put(8.4,1.15){$(a,e,b)$}
\end{picture}}
\caption{A fold of type IA}\label{ff1}
\end{figure}
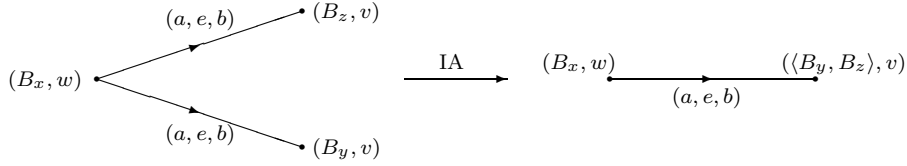

\medskip\noindent{\bf Fold of type IIIA}
  Suppose that $f_1$ and $f_2$ are both non-loop edges such that
  $y=\omega(f_1)=\omega(f_2)=z$. Furthermore $a=(f_1)_{\alpha}=(f_2)_{\alpha}$, $b=(f_1)_{\omega}$ and $b'=(f_2)_{\omega}$. Furthermore $e=]f_1]=[f_2]$, $a=(f_1)_{\alpha}=(f_2)_{\alpha}$, $b=(f_1)_{\omega}$ and $b'=(f_2)_{\omega}$.

  We identify the edges edges $f_1$ and $f_2$ of $\mathcal B'$ into a
  single edge $f$ with edge group $\langle B_{f_1},B_{f_2}\rangle$ and label $(a,e,b)$. We set the label of $\omega(f)$
  to be
\[
(\langle B_y, b^{-1}b'\rangle, v).\]

We call this operation a {\em fold of type}
IIIA.

\begin{figure}[here]
  \centerline{ \footnotesize \setlength{\unitlength}{.8cm}
\begin{picture}(10,3)
\bezier{200}(-2,1.5)(-.33,2.5)(.5,2.5)
\bezier{200}(-2,1.5)(-.33,.5)(.5,0.5)
\bezier{200}(.5,2.5)(1.33,2.5)(3,1.5)
\bezier{200}(.5,0.5)(1.33,.5)(3,1.5)
\put(4.5,1.5){\vector(1,0){1.5}}
\put(.5,2.5){\vector(1,0){.01}}
\put(.5,.5){\vector(1,0){.01}}
\put(9.5,1.5){\vector(1,0){.01}}
\put(7.5,1.5){\line(3,0){4}}
\put(-2,1.5){\circle*{.07}}
\put(3,1.5){\circle*{.07}}
\put(7.5,1.5){\circle*{.07}}
\put(11.5,1.5){\circle*{.07}}
\put(4.8,1.65){IIIA}
\put(-.3,2.65){$(a,e,b')$}
\put(-.3,0){$(a,e,b)$}
\put(6.5,1.05){$(B_x,w)$}
\put(10.3,1.1){$(\langle B_y,b^{-1}b'\rangle,v)$}
\put(8.7,1.75){$(a,e,b)$}
\put(-2.8,1){$(B_x,w)$}
\put(3.1,1.15){$(B_y,v)$}
\end{picture}}
\caption{A fold of type IIIA}\label{ff4}
\end{figure}
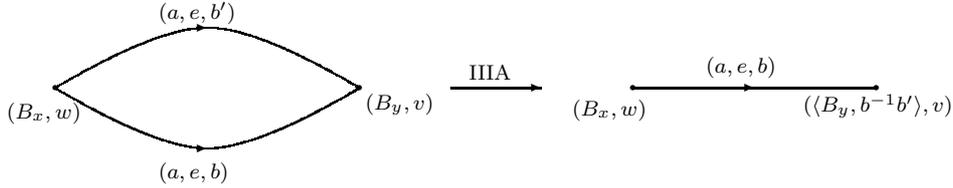

\medskip\noindent{\bf Fold of type IIA: }
  Let $\mathcal B$ be an $\mathbb A$-graph.  Suppose that $x\neq y$, i.e. that $f$ is a non-loop
  edge of $\mathcal B$ with the label $(a, e, b)$ and the edge group $B_f$. Suppose futher that $g\in A_{[f]}$ with $a\alpha_{e}(g)a^{-1}\in B_{\alpha(f)}$.

Let $\mathcal B'$ be the $\mathbb A$-graph obtained from $\mathcal
  B$ by replacing the the edge group $B_f$ by the group $\langle B_f,g\rangle$ and replacing the vertex group $B_y$ by $\langle B_y,b^{-1} \omega_e(g) b\rangle$.

We say that \emph{$\mathcal B'$ is obtained from $\mathcal B$ by a fold of type IIA}.

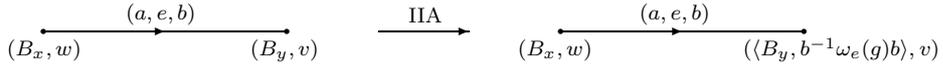
\begin{figure}[here]
  \centerline{ \footnotesize \setlength{\unitlength}{.8cm}
\begin{picture}(11,1)
\put(4.5,.5){\vector(1,0){1.5}}
\put(5,.65){IIA}
\put(-1,0.5){\line(3,0){4}}
\put(7.5,0.5){\line(3,0){4}}
\put(-1,.5){\circle*{.07}}
\put(3,.5){\circle*{.07}}
\put(7.5,.5){\circle*{.07}}
\put(11.5,.5){\circle*{.07}}
\put(-1.6,.1){$(B_x,w)$}
\put(6.8,.1){$(B_x,w)$}
\put(2.4,.1){$(B_y,v)$}
\put(10.5,.1){$(\langle B_y, b^{-1} \omega_e(g) b\rangle,v)$}
\put(.35,.7){$(a,e,b)$}
\put(8.8,.7){$(a,e,b)$}
\put(9.5,.5){\vector(1,0){.01}}
\put(1,.5){\vector(1,0){.01}}
\end{picture}}
\caption{A fold of type IIA with $B_f$ being replaced by $\langle B_f,g\rangle$}\label{ff5}
\end{figure}

Note that whenever $\mathcal B'$ is obtained from $\mathcal B$ by a fold as above then the fold induces a map $B\to B'$ of the underlying graphs, in the case of a fold of type IIA this is an isomorphism. We will denote this induced map by $\pi$.
\section{Proof of the theorem}\label{proof}

In the following we assume that $\mathbb A$ is a weakly reduced $(k,C)$-acylindrical graph of groups. Let $n$ be the rank of $\pi_1(\mathbb A)$. Thus we have to prove that $$\#EA\le (2k+1)\cdot C\cdot (n-1).$$ Note first that  we can assume that some vertex group of $\mathbb A$ is non-trivial. Indeed otherwise the we are in the graph setting and it is well known that a graph without valence 2 or 1 vertices and fundamental group of rank $n$ has at most $3n-3$ edges.

\medskip We study triples $(\mathcal B,\Gamma ,\mathcal E)$ where $\mathcal B$ is an $\mathbb A$-graph with associated graph of groups $\mathbb B$, $\Gamma $ is not necessarily connected subgraph of $B$ and $\mathcal E$ is a collection of edges of $B$ that is disjoint from $E\Gamma$ such that the following hold:

\begin{enumerate}
\item The complement of $\Gamma \cup \mathcal E$ in $B$ is a forest $\mathcal T$.
\item Every vertex of $B$ is either a vertex of $\Gamma $ or a vertex of some $T\in\mathcal T$.
\item Every $T\in\mathcal T$ meets $\Gamma $ in a single vertex $v_T$.
\item For every $T\in\mathcal T$ the vertex ${v_T}$ carries the fundamental group of $\mathbb T$, the subgraph of groups of $\mathbb B$ corresponding to $T$. Thus the inclusion $B_{v_T}\to \pi_1(\mathbb T,v_T)$ is surjective.
\end{enumerate}

We call such a triple $(\mathcal B,\Gamma ,\mathcal E)$ a {\em decorated} $\mathbb A$-graph. 

\medskip We will depict a decorated $\mathbb A$-graph such that the edges of $\Gamma $ are fat lines, the edges that belong to the trees $\mathcal T$ are thin lines that are oriented such that they point in $T$ away from $v_T$ if the edge lies in $T\in\mathcal T$. The remaining edges $\mathcal E$ are dotted lines. Components of $C$ that are single vertices are depicted by a fat dot.

\begin{figure}[here]
\footnotesize
\setlength{\unitlength}{.8cm}
\begin{center}
\begin{picture}(14,5)
%\bezier{15}(0,1.5)(.75,1.75)(1.5,2)
%\bezier{15}(0,.5)(.75,.25)(1.5,0)
\bezier{15}(9.5,0)(8.75,0)(8,0)
\bezier{15}(9.5,2)(8.75,2)(8,2)
\bezier{15}(12.5,1)(13.25,1.5)(14,2)
\bezier{25}(5,3)(6,3.5)(7,4)
\bezier{25}(5,5)(6,4.5)(7,4)
\bezier{20}(9.5,4)(9.5,3)(9.5,2)
%\put(1.5,0){\line(1,0){1.5}}
\bezier{15}(1.5,0)(2.25,0)(3,0)
\put(3,2){\line(0,1){2}}\put(3,3){\vector(0,1){.01}}
\put(3,4){\line(2,1){2}}\put(4,4.5){\vector(2,1){.01}}
\put(3,4){\line(2,-1){2}}\put(4,3.5){\vector(2,-1){.01}}
\put(6.5,1){\line(3,-2){1.5}}\put(7.25,.5){\vector(3,-2){.01}}
\put(9.5,2){\line(1,0){1.5}}\put(10.25,2){\vector(-1,0){.01}}
\put(12.5,1){\line(3,-2){1.5}}\put(13.25,.5){\vector(3,-2){.01}}
\linethickness{.6mm}
\bezier{100}(3,2)(2.25,2)(1.5,2)
\bezier{100}(3,2)(3,1)(3,0)
\bezier{100}(4.5,1)(3.75,.5)(3,0)
\bezier{100}(4.5,1)(3.75,1.5)(3,2)
\bezier{100}(4.5,1)(5.5,1)(6.5,1)
\bezier{100}(8,2)(7.25,1.5)(6.5,1)
\bezier{100}(9.5,0)(10.25,0)(11,0)
\bezier{100}(12.5,1)(11.75,.5)(11,0)
\bezier{100}(11,2)(11,1)(11,0)
\bezier{100}(11,2)(11.75,1.5)(12.5,1)
\bezier{100}(7,4)(8.25,4)(9.5,4)
\put(14,2){\circle*{.2}}
\put(1.5,0){\circle*{.2}}
\footnotesize
\end{picture}
\end{center}
\caption{A decorated $\mathbb A$-graph}
\end{figure}
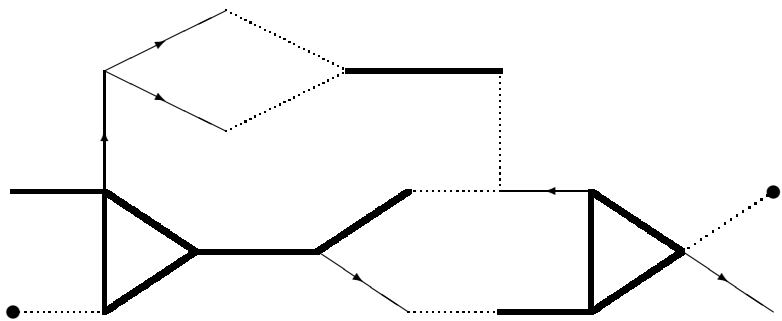

Suppose now that $(\mathcal B, \Gamma , \mathcal E)$ is a decorated $\mathbb A$-graph. Then we can associate to any vertex $v\in VB-V\Gamma $ a unique reduced path $\gamma_v$ such that $\gamma_v\subset T$ for some $T\in\mathcal T$, that $\alpha(\gamma_v)=v_T$ and $\omega(\gamma_v)=v$.

As $v_T$ carries the fundamental group of $\mathbb T$ it follows that for any $\gamma_v=e_1,\ldots ,e_k$ and $i=1,\ldots ,k$ the boundary monomorphism $\omega_{e_i}:B_{e_i}\to B_{\omega(e_i)}
$ is surjective. For any vertex of $v\in VB-V\Gamma $ will refer to the last edge of the path $\gamma_v$ by $e_v$. This implies in particular that the boundary monomorphism $\omega_{e_v}:B_{e_v}\to B_v$ is surjective for such vertices.

\medskip For any edge group $H$ we denote the number of times it can be replaced with a proper overgroup without yielding a group of order more than $C$ by $r(H)$. As replacing a finite group with a proper overgroup at least doubles the order it follows that for any finite group $H$ we have $r(H)\le\log_2\left(\frac{C}{|H|}\right)$. We further put $p(B)=2^{r(B)}$. Note that $p(B)=1$ if $B$ is of order more than $\frac{C}{2}$.

\medskip To any decorated $\mathbb A$-graph  $\tilde{\mathcal{B}}=(\mathcal B,\Gamma ,\mathcal E)$ we associate the complexity $$c(\tilde{\mathcal B})=\#[E\Gamma ]+(2k+1)\sum\limits_{e\in \mathcal E}p(B_e)$$ where $\#[E\Gamma ]$ is the number of edges of $A$ that lie in the image of $E\Gamma $ under the morphism $[\,.\,]:B\to A$.

%\medskip In the following we will refer to decorated $\mathbb A$-graphs which are triples of type $\tilde{\mathcal B}=(\mathcal B,\mathcal C,\{\gamma_v\,|\,v\in VB\})$ such that $\tilde{\mathbb B}=(\mathbb B,\mathcal C,\{\gamma_v\,|\,v\in VB\})$ is a decorated graph of groups where $\mathbb B$ is the graph of groups associated to $\mathcal B$. We put $c(\tilde{\mathcal B})=c(\tilde{\mathbb B})$, $Y(\tilde{\mathcal B})=Y(\tilde{\mathbb B})$, etc.

\medskip We denote by $\mathcal A$ is the trivial $\mathbb A$-graph, i.e. the $\mathbb A$-graph $\mathcal B$ whose graph morphism $[\,.\,]$ is an automorphism and where $B_v=A_{[v]}$ and $B_e=A_{[e]}$ for all $v\in VB$ and $e\in EA$.  Note that any two such $\mathbb A$-graphs are related by auxiliary moves. We first observe that a bound on the complexity of some decoration of $\mathcal A$ provides a bound on the complexity  of $\mathbb A$.

\begin{Lemma}\label{reduction} To prove Theorem~\ref{main} suffices to show that there exists a decoration $\tilde{\mathcal A}=(\mathcal A,\Gamma ,\mathcal E)$ of $\mathcal A$ such that $$c(\tilde{\mathcal A})\le (2k+1)C(n-1).$$\end{Lemma}

\noindent{\em Proof } The proof is by induction on the cardinality of $\mathcal E$ however  the claim that we prove by induction is slightly more complicated.  

Note first that that for any graph of groups $\mathbb B$ there exits a graph of groups $\mathbb B'$ obtained by ignoring vertices for which both boundary monomorphisms are surjective, i.e. by performing the inverse operation to a subdivision. For any graph of groups $\mathbb B$ we denote by $c_r(\mathbb B)$ the number of edges of $\mathbb B'$. Thus $\#EB=c_r(\mathbb B)$ iff $\mathbb B$ is weakly reduced, in particular we have $c_r(\mathbb A)=\#EA$. We call $c_r(\mathbb B)$ the {\em weakly reduced complexity} of $\mathbb B$.

Let $\mathcal E=\{e_1,\ldots ,e_l\}$. We define a filtration $X_0,X_1,\ldots ,X_l$ of $A$ such that $X_0=\Gamma $ and that $$X_i=X_{i-1}\cup \{e_i\}\cup \gamma_{\alpha(e_i)}\cup\gamma_{\omega(e_i)}$$ for $0\le i\le l$. Note that the definition of the decoration and the minimality of $\mathbb A$ imply that $X_l=A$. Denote the not necessarily connected subgraph of groups of $\mathbb A$ corresponding to $X_i$ by $\mathbb X_i$. We show by induction on $i$ that $$c_r(\mathbb X_i)\le \#[E\Gamma ]+(2k+1)\sum\limits_{j=1}^ip(A_{e_j}).$$ For $i=l$ this implies the assertion of the lemma as $$\#EA=c_r(\mathbb A)=c_r(\mathbb X_l)\le \#[E\Gamma ]+(2k+1)\sum\limits_{j=1}^lp(A_{e_j})=c(\tilde{\mathcal A})\le (2k+1)C(n-1).$$ The case $j=0$ is trivial as $\#[E\Gamma]=\#E\Gamma$. Recall that $X_i=X_{i-1}\cup \{e_i\}\cup \gamma_{\alpha(e_i)}\cup\gamma_{\omega(e_i)}$. Let $L_i= \{e_i\}\cup \gamma_{\alpha(e_i)}\cup\gamma_{\omega(e_i)}$ and $\mathbb L_i$ be the corresponding subgraph of groups of $\mathbb A$. 

Note that $c_r(\mathbb T_i)\le 2k+2p(A_{e_i})-1$. Indeed the first and the last $k$ edges of the path $\gamma_{\alpha(e_i)}\cup e_i\cup\gamma_{\omega(e_i)}^{-1}$ can have arbitrary edge groups but all edges in between have to have order at most $C$ due to the acylindricty assumption. If there any of the later edges then the order of the edge groups must decrease to  $|A_{e_i}|$ and then increase again. Thus there are at most $2p(A_{e_i})-1$ such edges. This implies the claim.

It is now easily verified that the inductive step follows. Note that when gluing $\gamma_{\alpha(e_i)}\cup e_i\cup\gamma_{\omega(e_i)}^{-1}$ to $X_{i-1}$ one potentially has to subdivide an edge first and therefore increase the reduced complexity. This however only happens if some other edge of $\gamma_{\alpha(e_i)}\cup e_i\cup\gamma_{\omega(e_i)}^{-1}$ lies in $X_{i-1}$ already as we assume that $k\ge 1$. Thus this does not affect the argument.\hfill$\Box$

\medskip The strategy of the proof is to first show that there exists a decorated $\mathbb A$-graph $\tilde{\mathcal B}$ of complexity $(2k+1)C(n-1)$ and then show that it can be folded onto a decoration of $\mathcal A$ without increasing the complexity. The assertion of the main theorem then follows from Lemma~\ref{reduction}. The first part is easy:

\begin{Lemma}\label{freeresolution} There exists a decorated $\mathbb A$-graph $\tilde{\mathcal B}_0$ that represents $\pi_1(\mathbb A)$ such that $c(\tilde{\mathcal B}_0)\le (2k+1)C(n-1).$
\end{Lemma}

\noindent{\em Proof} Choose an arbitrary minimal generating tuple $S=(g_1,\ldots, g_n)$ of $\pi_1(\mathbb A)$. As we are assuming that some vertex group of $\mathbb A$ is non-trivial it follows that after Nielsen-equivalence we can assume that $g_1$ is an elliptic element and that all other elements are hyperbolic.

We choose $\mathcal B_0$ to be the  $S$-wedge with $S$ as above.  Thus $\mathbb B_0$ is a graph of groups with a single cyclic vertex group, all other vertex and edge groups trivial and $(n-1)$ edges outside a maximal subtree. We decorate $\mathcal B_0$ by putting $\Gamma $ to be  the subgraph consisting of the single vertex with cyclic vertex group and choosing $\mathcal E$ to be the set of edges outside some maximal subtree. The complexity of the decorated $\mathbb A$-graph is cleary $0+(2k+1)(n-1)p(1)$ as there are $n-1$ edge outside a maximal tree and all edge groups are trivial. As $p(1)\le C$ this implies the assertion.\hfill$\Box$

\medskip The aim is to transform $\mathcal B_0$ into a folded $\mathbb A$-graph. This cannot always be done by finite folding sequences, in some situations we will therefore have to apply infinitely many folds at once.

%\medskip In the following we will see that we can fold decorated $\mathbb A$-graphs without increasing the complexity. At the end of the folding sequence we will have a folded decorated $\mathbb A$-graph that does not have larger complexity then the decorated $\mathbb A$-graph from Lemma~\ref{freeresolution}, the following Lemma shows that this proves our theorem.

\medskip In the following we denote the complement of $\mathcal E$ in $B$ by $Y=Y(\tilde{\mathcal B})$. % and the corresponding (not necessarily connected) subgraph of groups by $\mathbb Y=\mathbb Y(\tilde{\mathcal B})$. 
We call a decorated $\mathbb A$-graph $\tilde{\mathcal{B}}$ {\em tame} if the sub-$\mathbb A$-graphs whose underlying graphs are the components of $Y(\tilde{\mathcal{B}})$ are folded. 

Note that we can replace any decorated $\mathbb A$-graph $\tilde{\mathcal{B}}=(\mathcal B,\Gamma ,\mathcal E)$ by a tame decorated $\mathbb A$-graph $\tilde{\mathcal {B}}'=(\mathcal B',\Gamma ',\mathcal E')$ by applying all possible folding sequences to the components of $Y(\tilde{\mathcal B})$. The graph $\Gamma '$ and the set $\mathcal E'$ of the  decoration of $\mathcal B'$ are the images of $\Gamma $ and $\mathcal E$ under the folds. We call $\tilde{\mathcal {B}}'$ the {\em taming} of $\tilde{\mathcal {B}}$.

As these folds do not involve any of the edges of $\mathcal E$ and also preserves their edge groups and these folds further do not increased the number of edges of the graph $\Gamma $ we have the following:

\begin{Lemma}\label{lem:tame} Let $\tilde{\mathcal{B}}$ be a decorated $\mathbb A$-graph and  $\tilde{\mathcal {B}}'$ its taming. Then $c(\tilde{\mathcal {B}}')\le c(\tilde{\mathcal{B}})$.
\end{Lemma}

We will now show that we can assume that all edge groups of edges of $\mathcal E$ are of order at most $C$.

\medskip Let $\tilde{\mathcal B}=(\mathcal B,\Gamma ,\mathcal E)$ be a decorated $\mathbb A$-graph and $e\in \mathcal E$. Let $v=\alpha(e)$ and $w=\omega(e)$. We then associate to $\tilde{\mathcal B}$  the new decorated $\mathbb A$-graph $\tilde{\mathcal B}(e)=(\mathcal B,\Gamma ',\mathcal E')$ as follows (note that only the decoration changes):

\begin{enumerate}
\item $\Gamma '=\Gamma \cup e \cup \gamma_{\alpha(e)}\cup\gamma_{\omega(e)}$.
\item $\mathcal E'=\mathcal E-e$.%\gamma'_v$ is the suffix of $\gamma_v$ that meets no edge of $C'$.
\end{enumerate} 

%This definition implies in particular that $B-Y(\tilde{\mathcal B})=\left(B-Y(\tilde{\mathcal B}(e))\right)\cup e$.

\begin{Lemma}\label{lem:amalgamate} Let $\tilde{\mathcal B}$ be a tame decorated $\mathbb A$-graph and $e\in EB-EY$ such that $|B_e|>C$. Then $$c(\tilde{\mathcal B}(e))\le c(\tilde{\mathcal B}).$$
\end{Lemma}

\noindent{\em Proof} As $|B_e|>C$ it follows that $p(B_e)=1$. Thus when going from $\tilde{\mathcal B}$ to $\tilde{\mathcal B}(e)$ the second summand of the complexity decreases by $2k+1$. On the other hand we have $\#[E\Gamma ']\le \#[E\Gamma ]+2k+1$ as both $\gamma_{\alpha(v)}$ and $\gamma_{\omega(v)}$ lift to segments in the Bass-Serre tree that have stabilizer of order more than $C$ which implies that they are of projective length at most $k$ by assumption. Thus the assertion follows. \hfill$\Box$

\medskip We next see that the complexity is also well-behaved under folds. This is the main step in the proof of Theorem~\ref{main}.

\begin{Lemma}\label{lem:fold} Let $\tilde{\mathcal B}=(\mathcal B,\Gamma ,\mathcal E)$ be a tame decorated $\mathbb A$-graph such that all edges of $\mathcal E$ have order at most $C$. Suppose that $\mathcal B'$ is obtained from $\mathcal B$ by an elementary fold where we only apply folds of type I and III if no fold of type II can be applied.

Then there exists a decoration  $\tilde{\mathcal B}'=(\mathcal B',\Gamma ',\mathcal E')$ of $\mathcal B'$ such that $$c(\tilde{\mathcal B}')\le c(\tilde{\mathcal B}).$$
\end{Lemma}

\noindent{\em Proof} We deal with the different types of folds. We only discuss the A-type folds IA, IIA, IIIA for the $B$-type folds the arguments are very similar. It should also be noted that $A$ can have at most $n-1$ loop edges, thus we could just collapse those first and prove a slightly weaker bound on the complexity as there would be only folds of A-type then. Each type of fold has a couple of subcases depending on the local structure of the decoration. 

Note that we can assume that at least one edge affected by the fold lies in $\mathcal E$ because of the tameness assumption. In the case of folds of type IA and IIIA we can further assume that the two edges $e_1$ and $e_2$ have the same edge group which is also the edge group of the edge $\pi(e_1)=\pi(e_2)$ as otherwise a fold of type IIA would be possible. Note that we will in all cases describe the new decoration and we will leave the trivial observation that it is a decoration indeed to the reader.

\smallskip\noindent{\bf Folds of type IIA} The fold affects an edge $e\in \mathcal E$ with $|B_e|\le C$. Let $x=\alpha(e)$ and $y=\omega(e)$. Thus the fold adds an element $g\in B_x$ to $B_e$ and $B_y$. The following three cases clearly cover all possibilities.

\smallskip\noindent Case 1: $y\in V\Gamma $. In this case we keep the decoration. The complexity clearly does not increase as $B_e$ is replaced by a proper overgroup which does not increase the contribution of $e$ to the complexity. Note that the complexity does in fact decrease unless $|B_e|>\frac{1}{2}C$.

\smallskip\noindent Case 2: $y\notin V\Gamma $ and $\omega_{e}:B_{e}\to B_y$ is surjective. This implies in particular that $|B_e|=|B_y|=|B_{e_y}|$. The decoration of the new $\mathbb A$-graph is given by $\Gamma'=\Gamma$ and $\mathcal E'=(\mathcal E-e)\cup e_y$. As $|B_e|=|B_{e_y}|$ this does not change the complexity.

 
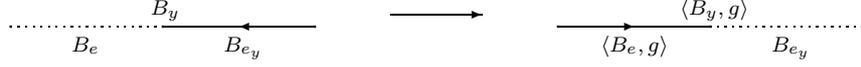
\begin{figure}[here]
\footnotesize
\setlength{\unitlength}{.8cm}
\begin{center}
\begin{picture}(14,1)
\bezier{20}(0,.3)(1.25,.3)(2.5,0.3)
\bezier{20}(14,.3)(12.75,.3)(11.5,0.3)
\put(6.25,.5){\vector(1,0){1.5}}
\put(2.5,.3){\line(1,0){2.5}}
\put(9,.3){\line(1,0){2.5}}
\put(10.25,.3){\vector(1,0){.01}}
\put(3.75,.3){\vector(-1,0){.01}}
\linethickness{.6mm}
\put(0,.3){\circle*{.05}}
\put(2.5,.3){\circle*{.05}}
\put(5,.3){\circle*{.05}}
\put(14,.3){\circle*{.05}}
\put(9,.3){\circle*{.05}}
\put(11.5,.3){\circle*{.05}}
\put(1,-.1){$B_e$}
\put(3.5,-.1){$B_{e_y}$}
\put(2.3,.5){$B_y$}
\put(9.7,-.1){$\langle B_e,g\rangle$}
\put(12.5,-.1){$B_{e_y}$}
\put(11,.5){$\langle B_y,g\rangle$}
\footnotesize
\end{picture}
\end{center}
\end{figure}

\smallskip\noindent Case 3: $y\notin V\Gamma $ and $\omega_{e}:B_{e}\to B_y$ is not surjective.

 If $|B_e|>\frac{1}{2}C$ then we  put $\tilde{\mathcal B}'=\tilde{\mathcal B}(e)$. The same argument as in the proof of Lemma~\ref{lem:amalgamate} shows that the complexity does not increase. 
 
 If $|B_e|\le \frac{1}{2}C$ then we put $\Gamma'=\Gamma\cup y $ which creates a new component consisting of a single vertex and put $\mathcal E'=\mathcal E\cup e_y$. As we are not adding an edge to $\Gamma $ this does not change $\#[E\Gamma ]$. In the sum part of the complexity the summand $(2k+1)\cdot p(B_e)$ is replaced by $(2k+1)\cdot (p(\langle B_e,g\rangle)+p(B_{e_y})$. As both $\langle B_e,g\rangle$ and $B_{e_y}$ are proper overgroups of $B_e$ this does not increase the complexity either as passing to a proper overgroup at least halves $p(B_e)$ since $|B_e|\le \frac{1}{2}C$.
 
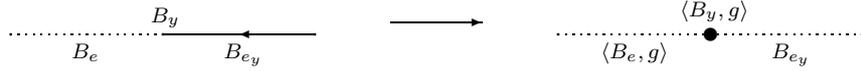
 \begin{figure}[here]
\footnotesize
\setlength{\unitlength}{.8cm}
\begin{center}
\begin{picture}(14,1)
\bezier{20}(0,.3)(1.25,.3)(2.5,0.3)
\bezier{40}(14,.3)(11.5,.3)(9,0.3)
\put(6.25,.5){\vector(1,0){1.5}}
\put(2.5,.3){\line(1,0){2.5}}
\put(3.75,.3){\vector(-1,0){.01}}
\linethickness{.6mm}
\put(0,.3){\circle*{.05}}
\put(2.5,.3){\circle*{.05}}
\put(5,.3){\circle*{.05}}
\put(14,.3){\circle*{.05}}
\put(9,.3){\circle*{.05}}
\put(11.5,.3){\circle*{.2}}
\put(1,-.1){$B_e$}
\put(3.5,-.1){$B_{e_y}$}
\put(2.3,.5){$B_y$}
\put(9.7,-.1){$\langle B_e,g\rangle$}
\put(12.5,-.1){$B_{e_y}$}
\put(11,.6){$\langle B_y,g\rangle$}
\footnotesize
\end{picture}
\end{center}
\end{figure}

\smallskip\noindent{\bf Folds of type IA} We will always assume that the decoration outside the portion of $B$ that is depicted is not changed. Note that we have $|B_{e_1}|=|B_{e_2}|=|B_{\pi(e_1)}|$. 

\smallskip 1)We first deal with the case that both $e_1$ and $e_2$ lie in $\mathcal E$.  There are five subcases, again they clearly deal with all possible situations. %In all four cases the edge $\obtained by identifying $e_1$ and $e_2$ lies in $B'-Y'$. 

1A) If both $y$ and $z$ are vertices of $\Gamma $ then we put $\Gamma'=\pi(\Gamma)$ and $\mathcal E'=\pi(\mathcal E)$. In this case the complexity clearly decreases.
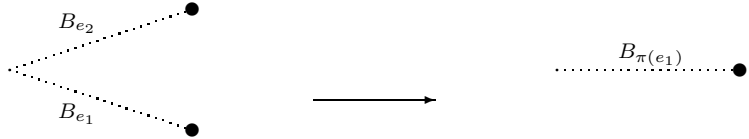
\begin{figure}[here]
\footnotesize
\setlength{\unitlength}{.8cm}
\begin{center}
\begin{picture}(14,2.3)
\bezier{25}(1,1)(2.5,.5)(4,0)
\bezier{25}(1,1)(2.5,1.5)(4,2)
\put(6,.5){\vector(1,0){2}}
\bezier{23}(10,1)(11.5,1)(13,1)
\linethickness{.6mm}
\put(4,0){\circle*{.2}}
\put(4,2){\circle*{.2}}
\put(13,1){\circle*{.2}}
\put(1,1){\circle*{.05}}
\put(10,1){\circle*{.05}}
\put(1.8,.2){$B_{e_1}$}
\put(1.8,1.7){$B_{e_2}$}
\put(11,1.2){$B_{\pi(e_1)}$}
\footnotesize
\end{picture}
\end{center}
\end{figure}

1B) $y$ is a vertex of $\Gamma $ and $z$ is not, the opposite case is analogous. Then $z$ has an associated edge $e_z$. We put $\Gamma'=\pi(\Gamma )$ and put $\mathcal E'=\pi(\mathcal E\cup e_z)$. As $|B_{e_z}|\ge |B_{e_1}|=|B_{e_2}|=|B_{\pi(e_1)}|$ this does not increase the complexity.

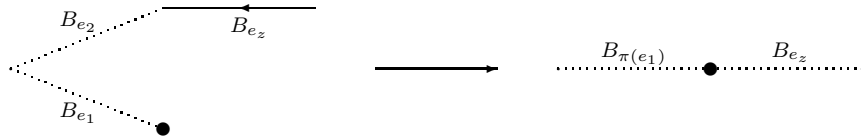
\begin{figure}[here]
\footnotesize
\setlength{\unitlength}{.8cm}
\begin{center}
\begin{picture}(14,2.3)
\bezier{22}(0,1)(1.25,.5)(2.5,0)
\bezier{22}(0,1)(1.25,1.5)(2.5,2)
\bezier{42}(9,1)(11.5,1)(14,1)
\put(6,1){\vector(1,0){2}}
\put(2.5,2){\line(1,0){2.5}}
\put(3.75,2){\vector(-1,0){.01}}
\linethickness{.6mm}
\put(2.5,0){\circle*{.2}}
\put(11.5,1){\circle*{.2}}
\put(0,1){\circle*{.05}}
\put(14,1){\circle*{.05}}
\put(9,1){\circle*{.05}}
\put(.8,.2){$B_{e_1}$}
\put(.8,1.7){$B_{e_2}$}
\put(3.6,1.6){$B_{e_z}$}
\put(9.7,1.2){$B_{\pi(e_1)}$}
\put(12.5,1.2){$B_{e_z}$}
\footnotesize
\end{picture}
\end{center}
\end{figure}

1C) Neither $y$ nor $z$ lie in $\Gamma $ and $\omega_{e_1}:B_{e_1}\to B_{\omega(e_1)}=B_y$ is surjective (the case that $\omega_{e_2}$ is surjective is analogous). Note that this implies that $|B_{e_1}|=|B_y|=|B_{e_y}|$. We put $\Gamma '=\pi(\Gamma )$ and $\mathcal E'=\pi(\mathcal E\cup e_y)$. As $|B_{e_1}|=|B_{e_2}|=|B_{e_y}|$ it follows that the complexity is unchanged.

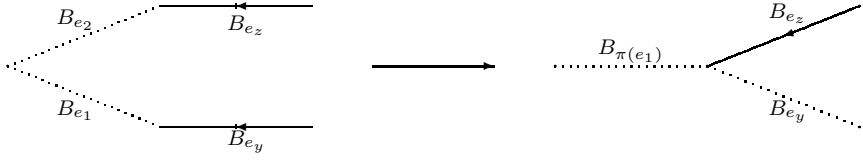
\begin{figure}[here]
\footnotesize
\setlength{\unitlength}{.8cm}
\begin{center}
\begin{picture}(14,2.3)
\bezier{22}(0,1)(1.25,.5)(2.5,0)
\bezier{22}(0,1)(1.25,1.5)(2.5,2)
\bezier{20}(9,1)(10.25,1)(11.5,1)
\bezier{200}(14,2)(12.75,1.5)(11.5,1)
\bezier{20}(14,0)(12.75,.5)(11.5,1)
\put(6,1){\vector(1,0){2}}
\put(2.5,2){\line(1,0){2.5}}
\put(2.5,0){\line(1,0){2.5}}
\linethickness{.6mm}
%\put(2.5,0){\circle*{.2}}
\put(11.5,1){\circle*{.05}}
\put(0,1){\circle*{.05}}
\put(14,2){\circle*{.05}}
\put(14,0){\circle*{.05}}
\put(9,1){\circle*{.05}}
\put(.8,.2){$B_{e_1}$}
\put(.8,1.7){$B_{e_2}$}
\put(3.6,1.6){$B_{e_z}$}
\put(3.6,-.3){$B_{e_y}$}
\put(9.7,1.2){$B_{\pi(e_1)}$}
\put(12.5,1.8){$B_{e_z}$}
\put(12.5,.2){$B_{e_y}$}
\put(3.75,0){\vector(-1,0){.01}}
\put(3.75,2){\vector(-1,0){.01}}
\put(12.75,1.5){\vector(-3,-1){.01}}
\footnotesize
\end{picture}
\end{center}
\end{figure}

1D) If neither $y$ nor $z$ lie in $\Gamma $, neiter $\omega_{e_1}:B_{e_1}\to B_{\omega(e_1)}=B_y$ nor $\omega_{e_2}:B_{e_2}\to B_z$ are surjective and $|B_{e_1}|=|B_{e_2}|\le \frac{1}{2}C$ then we put $\Gamma '=\pi(\Gamma \cup y)$ and $\mathcal E'=\pi(\mathcal E\cup \{e_y,e_z\})$. As both $B_{e_z}$ and $B_{e_y}$ are proper overgroups of $B_{e_1}$ it follows that the complexity does not increase.

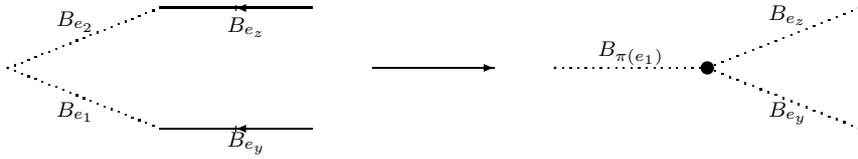
\begin{figure}[here]
\footnotesize
\setlength{\unitlength}{.8cm}
\begin{center}
\begin{picture}(14,2.3)
\bezier{22}(0,1)(1.25,.5)(2.5,0)
\bezier{22}(0,1)(1.25,1.5)(2.5,2)
\bezier{20}(9,1)(10.25,1)(11.5,1)
\bezier{20}(14,2)(12.75,1.5)(11.5,1)
\bezier{20}(14,0)(12.75,.5)(11.5,1)
\put(6,1){\vector(1,0){2}}
\put(2.5,2){\line(1,0){2.5}}
\put(2.5,0){\line(1,0){2.5}}
\linethickness{.6mm}
%\put(2.5,0){\circle*{.2}}
\put(11.5,1){\circle*{.2}}
\put(0,1){\circle*{.05}}
\put(14,2){\circle*{.05}}
\put(14,0){\circle*{.05}}
\put(9,1){\circle*{.05}}
\put(.8,.2){$B_{e_1}$}
\put(.8,1.7){$B_{e_2}$}
\put(3.6,1.6){$B_{e_z}$}
\put(3.6,-.3){$B_{e_y}$}
\put(9.7,1.2){$B_{\pi(e_1)}$}
\put(12.5,1.8){$B_{e_z}$}
\put(12.5,.2){$B_{e_y}$}
\put(3.75,0){\vector(-1,0){.01}}
\put(3.75,2){\vector(-1,0){.01}}
%\put(12.75,1.5){\vector(-3,-1){.01}}
\footnotesize
\end{picture}
\end{center}
\end{figure}

1E) If neither $y$ nor $z$ lie in $\Gamma $, neiter $\omega_{e_1}:B_{e_1}\to B_{\omega(e_1)}=B_y$ nor $\omega_{e_2}:B_{e_2}\to B_z$ are surjective and $|B_{e_1}|=|B_{e_2}|> \frac{1}{2}C$ then we put $\Gamma '=\pi(\Gamma \cup \gamma_y\cup\gamma_z)$ and $\mathcal E'=\pi(\mathcal E\cup \{e_y\})$. The same argument as in the proof of Lemma~\ref{lem:amalgamate} shows that the complexity does not increase.

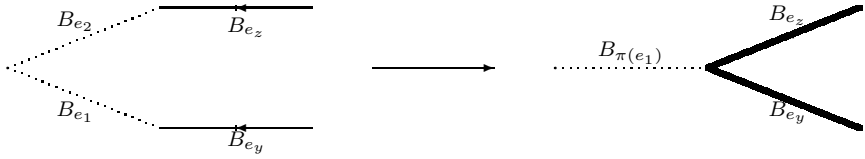
\begin{figure}[here]
\footnotesize
\setlength{\unitlength}{.8cm}
\begin{center}
\begin{picture}(14,2.3)
\bezier{22}(0,1)(1.25,.5)(2.5,0)
\bezier{22}(0,1)(1.25,1.5)(2.5,2)
\bezier{20}(9,1)(10.25,1)(11.5,1)
\put(6,1){\vector(1,0){2}}
\put(2.5,2){\line(1,0){2.5}}
\put(2.5,0){\line(1,0){2.5}}
\linethickness{.6mm}
%\put(2.5,0){\circle*{.2}}
%\put(11.5,1){\circle*{.2}}
\put(0,1){\circle*{.05}}
\put(14,2){\circle*{.05}}
\put(14,0){\circle*{.05}}
\put(9,1){\circle*{.05}}
\linethickness{.6mm}
\bezier{200}(14,2)(12.75,1.5)(11.5,1)
\bezier{200}(14,0)(12.75,.5)(11.5,1)
\put(.8,.2){$B_{e_1}$}
\put(.8,1.7){$B_{e_2}$}
\put(3.6,1.6){$B_{e_z}$}
\put(3.6,-.3){$B_{e_y}$}
\put(9.7,1.2){$B_{\pi(e_1)}$}
\put(12.5,1.8){$B_{e_z}$}
\put(12.5,.2){$B_{e_y}$}
\put(3.75,0){\vector(-1,0){.01}}
\put(3.75,2){\vector(-1,0){.01}}
\put(12.75,1.5){\vector(-3,-1){.01}}
\footnotesize
\end{picture}
\end{center}
\end{figure}

2) We next deal with the case that precisely one of the two edges $e_1$ and $e_2$, say $e_1$, lies in $\mathcal E$ and that $e_2\in E\Gamma$.

\smallskip 2A) If $y\in V\Gamma$ then we put $\Gamma'=\pi(\Gamma)$ and $\mathcal E'=\pi(\mathcal E-e_1)$, the complexity clearly decreases.

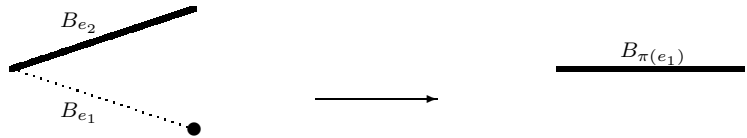
\begin{figure}[here]
\footnotesize
\setlength{\unitlength}{.8cm}
\begin{center}
\begin{picture}(14,2.3)
\put(6,.5){\vector(1,0){2}}
\bezier{23}(10,1)(11.5,1)(13,1)
\bezier{25}(1,1)(2.5,.5)(4,0)
\linethickness{.6mm}
\put(4,0){\circle*{.2}}
\put(4,2){\circle*{.05}}
\put(13,1){\circle*{.05}}
\put(1,1){\circle*{.05}}
\put(10,1){\circle*{.05}}
\put(1.8,.2){$B_{e_1}$}
\put(1.8,1.7){$B_{e_2}$}
\put(11,1.2){$B_{\pi(e_1)}$}
%\thicklines
\bezier{250}(1,1)(2.5,1.5)(4,2)
\bezier{230}(10,1)(11.5,1)(13,1)
\footnotesize
\end{picture}
\end{center}
\end{figure}

2B) If $y\notin\Gamma $ then we put put $\Gamma'=\pi(\Gamma )$ and put $\mathcal E'=\pi((\mathcal E-e_1)\cup e_y)$. As $|B_{e_y}|\ge |B_{e_1}|=|B_{e_2}|=|B_{\pi(e_1)}|$ this does not increase the complexity.

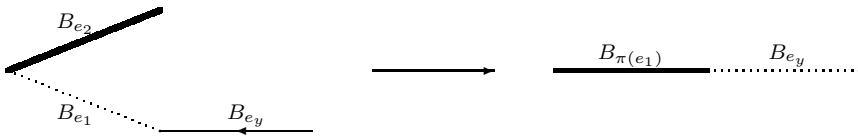
\begin{figure}[here]
\footnotesize
\setlength{\unitlength}{.8cm}
\begin{center}
\begin{picture}(14,2.3)
\bezier{21}(14,1)(12.75,1)(11.5,1)
\bezier{22}(0,1)(1.25,.5)(2.5,0)
\put(6,1){\vector(1,0){2}}
\put(2.5,0){\line(1,0){2.5}}
\put(3.75,0){\vector(-1,0){.01}}
\linethickness{.6mm}
\put(2.5,0){\circle*{.05}}
\put(11.5,1){\circle*{.05}}
\put(0,1){\circle*{.05}}
\put(14,1){\circle*{.05}}
\put(9,1){\circle*{.05}}
\put(.8,.2){$B_{e_1}$}
\put(.8,1.7){$B_{e_2}$}
\put(3.6,.2){$B_{e_y}$}
\put(9.7,1.2){$B_{\pi(e_1)}$}
\put(12.5,1.2){$B_{e_y}$}
\linethickness{.6mm}
\bezier{220}(0,1)(1.25,1.5)(2.5,2)
\bezier{210}(9,1)(10.25,1)(11.5,1)
\footnotesize
\end{picture}
\end{center}
\end{figure}

3) Precisely one of the two edges $e_1$ and $e_2$, say $e_1$, lies in $\mathcal E$ and that $e_2^{-1}$ is an edge of some path $\gamma_v$, i.e. the arrow on $e_2$ points towards $x$. In this case it turns out that we can argue precisely as in the case where $e_1,e_2\in\mathcal E$ except that this time the edge $\pi(e_1)=\pi(e_2)$ does not lie in $\mathcal E'$ but is an edge on some $\gamma_v$ pointing to $\pi(x)$.

\smallskip 4) We are left with the case that precisely one of the two edges $e_1$ and $e_2$, say $e_1$, lies in $\mathcal E$ and that $e_2$ is an edge of some path $\gamma_v$, i.e. the arrow on $e_2$ points towards $z$. Note that this implies that $z=\omega(e_2)\notin V\Gamma$. We distinguish the cases that $y\in V\Gamma$ and that $y\notin V\Gamma$.

\smallskip 4A) If $y\in V\Gamma$ then we put $\Gamma'=\pi(\Gamma)$ and $\mathcal E'=\pi(\mathcal E)$. The complexity is clearly unchanged.

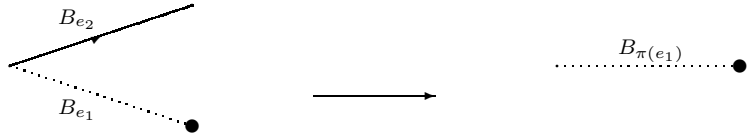
\begin{figure}[here]
\footnotesize
\setlength{\unitlength}{.8cm}
\begin{center}
\begin{picture}(14,2.3)
\put(6,.5){\vector(1,0){2}}
\put(2.5,1.5){\vector(2,1){.01}}
\bezier{25}(1,1)(2.5,.5)(4,0)
%\linethickness{.6mm}
\put(4,0){\circle*{.2}}
\put(4,2){\circle*{.05}}
\put(13,1){\circle*{.2}}
\put(1,1){\circle*{.05}}
\put(10,1){\circle*{.05}}
\put(1.8,.2){$B_{e_1}$}
\put(1.8,1.7){$B_{e_2}$}
\put(11,1.2){$B_{\pi(e_1)}$}
%\thicklines
\bezier{250}(1,1)(2.5,1.5)(4,2)
\bezier{23}(10,1)(11.5,1)(13,1)
\footnotesize
\end{picture}
\end{center}
\end{figure}

4B) If $y\notin\Gamma $ then we put put $\Gamma'=\pi(\Gamma )$ and put $\mathcal E'=\pi((\mathcal E-e_1)\cup e_y)$. As $|B_{e_y}|=|B_y|\ge |B_{e_1}|=|B_{e_2}|=|B_{\pi(e_1)}|$ this does not increase the complexity.

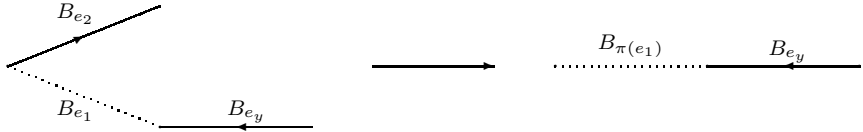
\begin{figure}[here]
\footnotesize
\setlength{\unitlength}{.8cm}
\begin{center}
\begin{picture}(14,2.3)
\put(1.25,1.5){\vector(2,1){.01}}
\bezier{21}(9,1)(10.25,1)(11.5,1)
\bezier{22}(0,1)(1.25,.5)(2.5,0)
\put(6,1){\vector(1,0){2}}
\put(2.5,0){\line(1,0){2.5}}
\put(3.75,0){\vector(-1,0){.01}}
\put(12.75,1){\vector(-1,0){.01}}
%\linethickness{.6mm}
\put(2.5,0){\circle*{.05}}
\put(11.5,1){\circle*{.05}}
\put(0,1){\circle*{.05}}
\put(14,1){\circle*{.05}}
\put(9,1){\circle*{.05}}
\put(.8,.2){$B_{e_1}$}
\put(.8,1.8){$B_{e_2}$}
\put(3.6,.2){$B_{e_y}$}
\put(9.7,1.3){$B_{\pi(e_1)}$}
\put(12.5,1.2){$B_{e_y}$}
%\linethickness{.6mm}
\bezier{220}(0,1)(1.25,1.5)(2.5,2)
\bezier{210}(14,1)(12.75,1)(11.5,1)
\footnotesize
\end{picture}
\end{center}
\end{figure}

\smallskip\noindent{\bf Folds of type IIIA} 

1) Again we deal first with the case that both $e_1$ and $e_2$ lie in $\mathcal E$. There are two cases to consider.

1A) If $y\in \Gamma $ then we put $\Gamma '=\pi(\Gamma )$ and $\mathcal E'=\pi(\mathcal E)$. Clearly the complexity decreases.

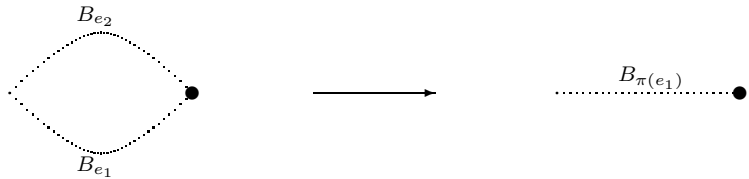
\begin{figure}[here]
\footnotesize
\setlength{\unitlength}{.8cm}
\begin{center}
\begin{picture}(14,2.5)
\bezier{20}(1,1)(2,2)(2.5,2)
\bezier{20}(1,1)(2,0)(2.5,0)
\bezier{20}(2.5,2)(3,2)(4,1)
\bezier{20}(2.5,0)(3,0)(4,1)
\bezier{28}(10,1)(11.5,1)(13,1)
\put(6,1){\vector(1,0){2}}
\linethickness{.6mm}
%\put(2.5,0){\circle*{.2}}
\put(13,1){\circle*{.2}}
\put(1,1){\circle*{.05}}
\put(4,1){\circle*{.2}}
\put(10,1){\circle*{.05}}
\put(2.1,-.3){$B_{e_1}$}
\put(2.1,2.2){$B_{e_2}$}
\put(11,1.2){$B_{\pi(e_1)}$}
\footnotesize
\end{picture}
\end{center}
\end{figure}

1B) If $y\notin \Gamma $ then we put $\Gamma '=\pi(\Gamma \cup y)$ and $\mathcal E'=\pi(\mathcal E\cup e_y)$. As $|B_{e_y}|\ge |B_{e_1}|$ it follows that the complexity does not increase.

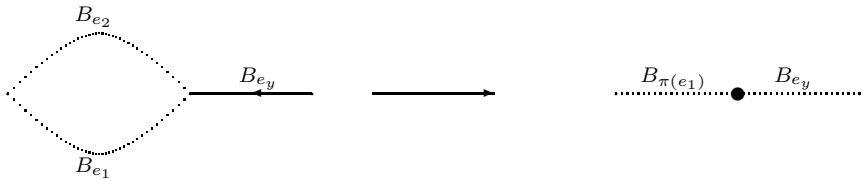
\begin{figure}[here]
\footnotesize
\setlength{\unitlength}{.8cm}
\begin{center}
\begin{picture}(14,2.5)
\bezier{20}(0,1)(1,2)(1.5,2)
\bezier{20}(0,1)(1,0)(1.5,0)
\bezier{20}(1.5,2)(2,2)(3,1)
\bezier{20}(1.5,0)(2,0)(3,1)
\bezier{44}(10,1)(12,1)(14,1)
\put(3,1){\line(1,0){2}}
\put(4,1){\vector(-1,0){.01}}
\put(6,1){\vector(1,0){2}}
\linethickness{.6mm}
%\put(2.5,0){\circle*{.2}}
\put(12,1){\circle*{.2}}
\put(0,1){\circle*{.05}}
\put(14,1){\circle*{.05}}
\put(3,1){\circle*{.05}}
\put(5,1){\circle*{.05}}
\put(10,1){\circle*{.05}}
\put(1.1,-.3){$B_{e_1}$}
\put(1.1,2.2){$B_{e_2}$}
\put(10.4,1.2){$B_{\pi(e_1)}$}
\put(3.8,1.2){$B_{e_y}$}
\put(12.6,1.2){$B_{e_y}$}
\footnotesize
\end{picture}
\end{center}
\end{figure}

2) We are left with the case that precisely one of the two edges $e_1$ and $e_2$ is in $\mathcal E$. Suppose that $e_2\notin\mathcal E$. 

\smallskip 2A) If $e_2\in E\Gamma$ then we put $\Gamma'=\pi(\Gamma)$ and $\mathcal E'=\pi(\mathcal E-e_1)$. The complexity clearly decreases.

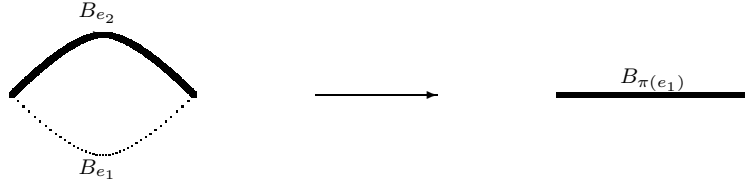
\begin{figure}[here]
\footnotesize
\setlength{\unitlength}{.8cm}
\begin{center}
\begin{picture}(14,2.5)
\bezier{20}(1,1)(2,0)(2.5,0)
\bezier{20}(2.5,0)(3,0)(4,1)
\put(6,1){\vector(1,0){2}}
\put(2.5,2){\vector(1,0){.01}}
\linethickness{.6mm}
%\put(2.5,0){\circle*{.2}}
%\put(13,1){\circle*{.05}}
\put(1,1){\circle*{.05}}
\put(4,1){\circle*{.05}}
%\put(10,1){\circle*{.05}}
\put(2.1,-.3){$B_{e_1}$}
\put(2.1,2.3){$B_{e_2}$}
\put(11,1.2){$B_{\pi(e_1)}$}
\linethickness{.6mm}
\bezier{200}(1,1)(2,2)(2.5,2)
\bezier{200}(2.5,2)(3,2)(4,1)
\bezier{280}(10,1)(11.5,1)(13,1)
\footnotesize
\end{picture}
\end{center}
\end{figure}

2B) If $y\in V\Gamma$ and $e_2^{-1}$ is an edge of some path $\gamma_v$, i.e. the arrow on $e_2$ points towards the vertex $x$ then we again put $\Gamma'=\pi(\Gamma)$ and $\mathcal E'=\pi(\mathcal E-e_1)$. Again the complexity decreases.

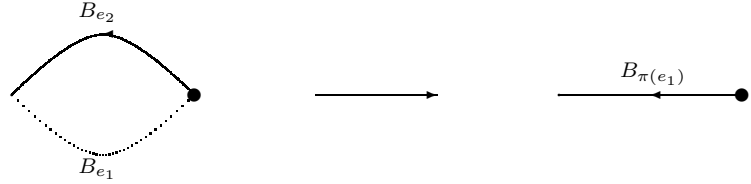
\begin{figure}[here]
\footnotesize
\setlength{\unitlength}{.8cm}
\begin{center}
\begin{picture}(14,2.5)
\bezier{200}(1,1)(2,2)(2.5,2)
\bezier{20}(1,1)(2,0)(2.5,0)
\bezier{200}(2.5,2)(3,2)(4,1)
\bezier{20}(2.5,0)(3,0)(4,1)
\bezier{280}(10,1)(11.5,1)(13,1)
\put(6,1){\vector(1,0){2}}
\put(11.5,1){\vector(-1,0){.01}}
\put(2.5,2){\vector(-1,0){.01}}
\linethickness{.6mm}
%\put(2.5,0){\circle*{.2}}
\put(13,1){\circle*{.2}}
\put(1,1){\circle*{.05}}
\put(4,1){\circle*{.2}}
\put(10,1){\circle*{.05}}
\put(2.1,-.3){$B_{e_1}$}
\put(2.1,2.3){$B_{e_2}$}
\put(11,1.3){$B_{\pi(e_1)}$}
\footnotesize
\end{picture}
\end{center}
\end{figure}

2C) If $y\notin V\Gamma$ and $e_2^{-1}$ is an edge of some path $\gamma_v$, i.e. the arrow on $e_2$ points towards the vertex $x$ then we again put $\Gamma'=\pi(\Gamma\cup y)$ and $\mathcal E'=\pi((\mathcal E\cup e_y)-e_1)$.  As $|B_{e_y}|\ge |B_e|$ it follows that the complexity does not increase.

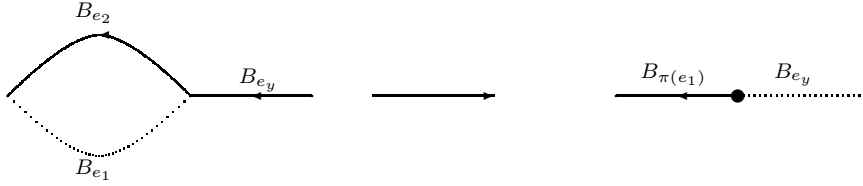
\begin{figure}[here]
\footnotesize
\setlength{\unitlength}{.8cm}
\begin{center}
\begin{picture}(14,2.5)
\bezier{200}(0,1)(1,2)(1.5,2)
\bezier{20}(0,1)(1,0)(1.5,0)
\bezier{200}(1.5,2)(2,2)(3,1)
\bezier{20}(1.5,0)(2,0)(3,1)
\bezier{22}(12,1)(13,1)(14,1)
\put(10,1){\line(1,0){2}}
\put(3,1){\line(1,0){2}}
\put(4,1){\vector(-1,0){.01}}
\put(6,1){\vector(1,0){2}}
\put(1.5,2){\vector(-1,0){.01}}
\put(11,1){\vector(-1,0){.01}}
\linethickness{.6mm}
%\put(2.5,0){\circle*{.2}}
\put(12,1){\circle*{.2}}
\put(0,1){\circle*{.05}}
\put(14,1){\circle*{.05}}
\put(3,1){\circle*{.05}}
\put(5,1){\circle*{.05}}
\put(10,1){\circle*{.05}}
\put(1.1,-.3){$B_{e_1}$}
\put(1.1,2.3){$B_{e_2}$}
\put(10.4,1.3){$B_{\pi(e_1)}$}
\put(3.8,1.2){$B_{e_y}$}
\put(12.6,1.3){$B_{e_y}$}
\footnotesize
\end{picture}
\end{center}
\end{figure}

2D) Thus we are left with the case that $e_2$ is an edge of some path $\gamma_v$, i.e. that the arrow on $e_2$ points towards $y$, in particular $y\notin V\Gamma$. In this case we put $\Gamma'=\pi(\Gamma\cup y)$ and $\mathcal E'=\pi(\mathcal E)$. The complexity is unchanged.

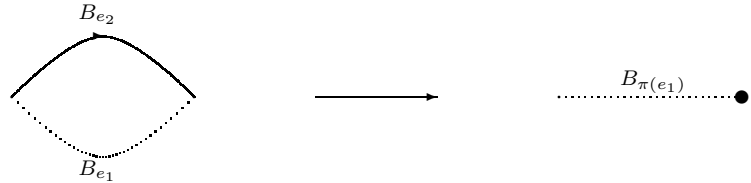
\begin{figure}[here]
\footnotesize
\setlength{\unitlength}{.8cm}
\begin{center}
\begin{picture}(14,2.5)
\bezier{200}(1,1)(2,2)(2.5,2)
\bezier{20}(1,1)(2,0)(2.5,0)
\bezier{200}(2.5,2)(3,2)(4,1)
\bezier{20}(2.5,0)(3,0)(4,1)
\bezier{28}(10,1)(11.5,1)(13,1)
\put(6,1){\vector(1,0){2}}
\put(2.5,2){\vector(1,0){.01}}
\linethickness{.6mm}
%\put(2.5,0){\circle*{.2}}
\put(13,1){\circle*{.2}}
\put(1,1){\circle*{.05}}
\put(4,1){\circle*{.05}}
\put(10,1){\circle*{.05}}
\put(2.1,-.3){$B_{e_1}$}
\put(2.1,2.3){$B_{e_2}$}
\put(11,1.2){$B_{\pi(e_1)}$}
\footnotesize
\end{picture}
\end{center}
\end{figure}\hfill$\Box$

\medskip We now have all necessary tools to conclude.

\medskip \noindent{\em Proof of Theorem~\ref{main}} Note first that there exists a finite sequence of decorated $\mathbb A$-graphs $$\tilde{\mathcal B}_0,\tilde{\mathcal B}_0',\tilde{\mathcal B}_1,\tilde{\mathcal B}_1',\tilde{\mathcal B}_2,\tilde{\mathcal B}_2',\ldots ,\tilde{\mathcal B}_l,\tilde{\mathcal B}_l'$$ such that the following hold:

\begin{enumerate}
\item $\tilde{\mathcal B}_0$ is as in the conlusion of Lemma~\ref{freeresolution}.
\item  $\tilde{\mathcal B}_l'$ is folded and surjective, i.e. $\mathbb B'_l\cong\mathbb A$. 
\item $\tilde{\mathcal B}'_{i}$ is the taming of $\tilde{\mathcal B}_i$ for all $i$.
\item $\tilde{\mathcal B}_i$ is obtained from $\tilde{\mathcal B}_{i-1}'$ by one of the following operations:
\begin{enumerate}
\item $\tilde{\mathcal B}_i=\tilde{\mathcal B}_{i-1}'(e)$ for some $e\in \mathcal E_{i-1}'$ with $|B_e|>C$.
\item $\tilde{\mathcal B}_i$ is obtained from $\tilde{\mathcal B}_{i-1}'$ by an elementary fold.
\end{enumerate}
\end{enumerate}

To see this start with $\tilde{\mathcal B}_0$ and apply step 4(a) followed by a taming as long as possible. This can only happen finitely many times. If not we apply step 4(b), i.e. apply an elementary fold again followed by a taming. There are only finitely many folds of type I and III as those decrease the complexity of the graph and folds of type II increase the order of an edge group that is of order less than $C$. This can also only happen finitely many times. Thus the process must stop with a folded decorated $\mathbb A$-graph.

\smallskip It follows from the Lemma~\ref{freeresolution}, \ref{lem:tame}, \ref{lem:amalgamate} and \ref{lem:fold}  that $$(2k+1)C(n-1)\ge c(\tilde{\mathcal B}_0)\ge c(\tilde{\mathcal B}_l)$$ and as $\mathbb B_l'$ is isomorphic to $\mathbb A$ this inequality implies the theorem by Lemma~\ref{reduction}.\hfill$\Box$

% with the property that for any vertex $v$ there exists a path $\gamma=e_1,\ldots ,e_k$ in $B$ that originates in some graph $\Gamma\in\mathcal C$ and ends at $v$ such that $\omega_{e_i}$ is surjective for i=1,\ldots ,k$.

% Suppose now that $\mathbb B$ is a graph of groups with underlying graph $B$. Suppose further 
%\section{Remarks}

\bibliographystyle{amsplain}

Department of Mathematics, Heriot-Watt University, Riccarton, Edinburgh EH14 4AS, Scotland, UK.

{\bf email: }R.Weidmann@ma.hw.ac.uk

\end{document}